\documentclass[11pt,a4paper]{amsart}
\usepackage[latin1]{inputenc}
\usepackage{amssymb}
\usepackage{amsmath}

\newtheorem{theorem}{Theorem}[section]
\newtheorem{rem}[theorem]{Remark}
\newtheorem{prop}[theorem]{Proposition}
\newtheorem{lemma}[theorem]{Lemma}
\newtheorem{cor}[theorem]{Corollary}
\newtheorem{defi}[theorem]{Definition}

\newcommand{\bprf}{{\it Proof.~}}
\newcommand{\eprf}{\hfill $\square$ \bigskip\par}
\newcommand{\erem}{\hfill $\square$}
\newcommand{\PP}{ \mathbb{P}}

\newcommand{\pdual}{(\mathbb{P}^n)^*}
\newcommand{\ch}{\mathcal{H}}

\input epsf
\input cyracc.def

\makeatletter
\def\blfootnote{\xdef\@thefnmark{}\@footnotetext}

\begin{document}


\title{A Geometrical approach to Gordan--Noether's and Franchetta's contributions to a question posed by Hesse.}
\author{Alice Garbagnati and Flavia Repetto}
\address{Alice Garbagnati, Dipartimento di Matematica, Universit\`a di Milano,
  via Saldini 50, I-20133 Milano, Italia}

\email{garba@mat.unimi.it}

\address{Flavia Repetto, Dipartimento di Matematica, Universit\`a di Milano,
  via Saldini 50, I-20133 Milano, Italia}

\email{flaviarepetto@yahoo.it}

\maketitle \markboth{ALICE GARBAGNATI AND FLAVIA REPETTO}{A
GEOMETRICAL APPROACH TO A QUESTION POSED BY HESSE.}

\begin{abstract}
Hesse claimed in \cite{Hesse1} (and later also in \cite{Hesse2}) that an irreducible
projective hypersurface in $\PP^n$ defined by an equation with vanishing hessian determinant
is necessarily a cone. Gordan and Noether proved in \cite{Gordan Noether} that this is true
for $n\leq 3$ and constructed counterexamples for every $n\geq 4$. Gordan and Noether
 and Franchetta  gave classification of hypersurfaces
in $\PP^4$ with vanishing hessian and which are not cones, see
\cite{Gordan Noether} and \cite{Franchetta}. Here we translate in
geometric terms Gordan and Noether approach, providing direct
geometrical proofs of these results.
\end{abstract}

\section{Introduction}
Let $f=f(x_0, \ldots, x_n)\in k[x_0, \ldots, x_n]$ be a non-zero
irreducible homogeneous polynomial  over an algebraically closed
field $k$ of characteristic zero. Then the hessian polynomial of
$f$ is the determinant of the matrix of the second derivatives:
 $$h_f:={\rm det} ([\partial^2
f/\partial x_i\partial x_j]_{i,j=0,\ldots ,n}).$$ Obviously when
the hypersurface $X=V(f)\subset\mathbb{P}^n$ is a cone (i.e. up
to a linear change of coordinates $f$ does not depend on all the
variables) then the hessian polynomial of $f$ is identically zero.
The converse is clearly true when ${\rm deg}(f)=2$.
Hesse claimed twice that the converse is true for
each degree of the polynomial $f$, i.e. he claimed that if the
hessian polynomial of a polynomial $f$ is identically zero then
the hypersurface $X=V(f)\subset\mathbb{P}^n$ is a cone (see
\cite{Hesse1}, \cite{Hesse2}).

The problem was reconsidered by Gordan and Noether (\cite{Gordan
Noether}) who proved that Hesse's claim is true when $n\leq 3$ but
false in general when $n\geq 4$. They constructed families of
counterexamples for every $n\geq 4$, which have been revisited
recently by Permutti in \cite{Permutti1}, \cite{Permutti2} and by
Lossen in \cite{Lossen}. Moreover, Gordan and Noether seem to have proved that their
families of examples  are the only possible counterexamples if $n=4$ but it is
rather difficult to indicate a precise reference for this result in their monumental paper.
Franchetta (\cite{Franchetta}) gave an independent classification
of hypersurfaces in $\mathbb{P}^4$ with vanishing hessian which
are not cones using more geometrical techniques. Other examples
were given by Perazzo, \cite{Perazzo}, who considered the case of
cubic hypersurfaces with vanishing hessian and obtained
the classification of these cubics in $\PP^4$, $\PP^5$ and
$\PP^6$. \\
Since the problem posed by Hesse has a  geometrical flavour, the
aim of this note is to translate in more geometric term Gordan and
Noether approach using some ideas and results contained in
\cite{Gordan Noether} and \cite{Lossen} and in the recent
\cite{CilibertoRussoSimis}. We also briefly describe the
counterexamples in projective spaces of dimension at least four
produced by Gordan and Noether, relating them  to works of
Franchetta and Permutti and we will give a short geometrical proof
of the characterization of hypersurfaces in $\mathbb{P}^4$  with
vanishing hessian which are not cones.

In the first Section we describe some background materials and we
consider a geometrical construction involving the dual variety of
a hypersurface. This construction allows us to reconsider the
Gordan and Noether's results and to describe them in a geometrical
context. In the second Section the Hesse's claim is proved in the
case of hypersurface of dimension at most 2. This proof is very
easy and it is based on the geometrical construction given in the
first Section. In the third Section the counterexamples by Gordan
and Noether and Franchetta are described, using also the results
of \cite{Permutti1}, \cite{Permutti2}, \cite{CilibertoRussoSimis}.
The last Section is dedicated to  hypersurfaces in $\PP^4$. We
describe the properties of  hypersurfaces in $\PP^4$ with
vanishing hessian and then we give another proof of  Franchetta's
classification of these hypersurfaces.\\
\medskip

{\bf Acknowledgements}.
We started our collaboration on this subject at Pragmatic 2006.
We would like to thank the organizers for the event and
  Professor Francesco Russo, who presented us the
problem and helped us during the preparation of this paper with
many corrections and suggestions.

\section{Background material.}

\subsection{The Polar map and the Hessian  of a projective hypersurface.}

Consider a non-constant homogeneous polynomial of degree $d\geq
1$, $f=f(x_0,\ldots,x_n)\in k[x_0, \ldots, x_n]$, in the $n+1$
variables $x_0, \ldots, x_n$ over an algebraically closed field
$k$ of characteristic zero. Denote by $f_i$  the partial
derivatives $\partial f/\partial x_i$, $i=0,\ldots ,n$.

\begin{defi}{\em Let $X=V(f)\subset\PP^n$ be the associated  hypersurface. We say that $X$
is a {\it cone} if, modulo projective transormations of  $\mathbb{P}^n$, the equation defining $X$ does not depend on all the  variables.

 Equivalently $X$ is a cone if and only if ${\rm Vert}(X)\neq\varnothing$. The {\it vertex of $X$}, ${\rm Vert}(X)$, is the set:
$${\rm Vert}(X):=\{ x\in X:\quad J(x,X)=X\},$$
where
 $$J(x,X)=\overline{\bigcup_{y\neq x, y\in X}<x,y>}\subset \mathbb{P}^n$$
 is {\it the join of $x$ and $X$}.
}\end{defi}
We recall that, if $X\subset\mathbb{P}^n$ is an (irreducible) subvariety of ${\rm dim}(X)=d$, then
$${\rm Vert}(X)=\bigcap_{x\in X}T_xX=\mathbb{P}^l\subset X,$$
with $l\geq -1$. (see e.g. \cite{Russo}, Proposition 1.2.6).

\begin{defi}{\em The {\it $($first$)$ polar map} associated to the hypersurface $X=V(f)\subset\mathbb{P}^n$ is the rational map
$\phi_f:\mathbb{P}^n\dashrightarrow\mathbb{P}^n$, defined by the partial derivatives of $f$:
$$\phi_f(p)= (f_0(p), \ldots, f_n(p)).$$
If $p\in X=V(f)$ is not singular, the polar map $\phi_f$ can be
interpreted as mapping the point $p\in X$  to its tangent
hyperplane $T_{X,p}$ (and, as such, the target of the map $\phi_f$
is $\mathbb{P}^{n*}$ ). Note that the base locus of $\phi_f$ is
the scheme ${\rm Sing}(X)=V(f_0, \ldots, f_n) \subset\PP^n$.
Denote by $Z(f)\subset\PP^{n*}$ the closure of the image of
$\mathbb{P}^n$ under the polar map $\phi_f$. The variety
$Z(f)\subset\mathbb{P}^{n*}$  is called the { \it polar image} of
$f$.}\end{defi}

\begin{defi}{\em We define   the {\it Hessian matrix} of the polynomial
$f$ to be  the $(n+1)\times(n+1)$ matrix:
$$\ch_f :=\left( \frac{\partial^2 f}{\partial x_i\partial x_j}\right) _{i,j=0,\ldots ,n}.$$
Its determinant $h_f:={\rm det}(\ch_f)\in
k[x_0, \ldots, x_n]$ is  the {\it Hessian polynomial} of
$f$.}\end{defi}

Note that the Jacobian matrix $J_{\phi_f}$
of  the (affine) polar map $\phi_f$ is exactly the Hessian matrix of $f$,
$\ch_f$.

We recall now the construction of the dual variety $X^*$ of an
algebraic reduced  variety $X\subset \mathbb{P}^n$.

Let ${\rm Sm}(X)$ denote the open non-empty subset of non singular points of a reduced variety $X\subset \mathbb{P}^n$.
Let
$$\mathcal{P}_X:=\overline{\{(x,H) : x\in {\rm Sm}(X), T_xX\subset H\}}\subset X\times(\mathbb{P}^n)^*$$ be the conormal variety of $X$, and consider the projections of $\mathcal{P}_X$ onto the factors:
$$
\begin{array}{ccccl}
&&\mathcal{P}_X&&\\
p_1&\swarrow&&\searrow&p_2\\
X&&&&X^*=\overline{p_2(\mathcal{P}_X)}\subset(\mathbb{P}^n)^*.\\
\end{array}
$$
The dual variety of $X$, $X^*$, is the scheme-theoretic image of
$\mathcal{P}_X$ in $(\mathbb{P}^n)^*$. In particular if $X$ is a
hypersurface of $\mathbb{P}^n$, then $X^*$ is the closure of the
set of hyperplanes tangent to $X$ at non--singular points. Observe
that  since  the Gauss map of $X$ associates to a non singular
point  $p\in X$ the point in $(\mathbb{P}^n)^{*}$ corresponding to
the hyperplane tangent to $X$ in $p$, we infer that (when $X$ is a
hypersurface) the closure of the image of the Gauss map of $X$ is
exactly the dual variety $X^*$.

Note also that the restriction of the polar map $\phi_f$ to
$V(f)\setminus{\rm Sing} (V(f))$  is the Gauss map of $X=V(f)$,
hence the closure of the  image of $X$ via $\phi_f$ is the dual
variety $X^*$ of $X$.

\subsection{Hypersurfaces with vanishing Hessian.}

We recall that $f_0, \ldots, f_n$ are \textit{algebraically
dependent} if there exists a polynomial $$\pi=\pi (y_0, \ldots,
y_n)\in k[y_0, \ldots, y_n]$$ such that $\pi(f_0,\ldots, f_n)=0$.
In particular they are linearly dependent if and only if there
exists such a $\pi$ of degree one.

Note first that the following  easy fact holds, recalling that the
Jacobian matrix of the affine polar map
$\widehat{\phi_f}:k^{n+1}\dasharrow k^{n+1}$ is the hessian matrix
$\mathcal H_f$.

\begin{prop}\label{prop: h=0 iff non dominant} Let $f\in k[x_0, \ldots, x_n]$ be an homogeneous polynomial, then the following are equivalent:
\begin{itemize}
\item $h_f\equiv 0$; \item $\phi_f$ is not a dominant map; \item
$Z(f)\subsetneq\PP^{n*}$; \item $f_0,\ldots, f_n$ are
algebraically dependent.
\end{itemize}
\end{prop}

We recall the following result from \cite{DimcaPapadima}, which
proves a conjecture stated in \cite{Dolg}.

\begin{theorem}\cite[Corollary 2]{DimcaPapadima}\label{polar degree}
The degree of the polar map $\phi_f$ depends only on ${\rm
Supp}(V(f))$ $($where the degree of  $\phi_f$ is zero if and only
if $\phi_f$ is not a dominant map$)$.
\end{theorem}

Note that by Proposition \ref{prop: h=0 iff non dominant} the
property of having vanishing Hessian is equivalent to the fact
that ${\rm dim}(Z(f))<n$, whence by Theorem \ref{polar degree} this
property depends only on the support of the hypersurface $X=V(f)$.

Since we are interested in hypersurfaces with vanishing Hessian,
from now on we shall assume that $X=V(f)$ is a reduced (and
irreducible) hypersurface.\\

The following result is due to Zak (see \cite{Zak}, Proposition
4.9).
\begin{prop}\label{prop: dimension of Z dual}
Let $X=V(f)\subset \mathbb{P}^n$ be a reduced hypersurface with
vanishing Hessian and let $Z(f)\subsetneq \mathbb{P}^{n*}$ denote
the polar image of $f$. Suppose $d\geq 2$, i.e. $\phi_f$ not
constant. Then:
$$Z(f)^*\subset {\rm Sing} (X).$$
In particular, ${\rm Sing(X)}\neq\emptyset$, ${\rm dim}
(Z(f)^*)\leq n-2$  and $X^*\subsetneq Z(f)$.
\end{prop}

In the sequel  we shall need the following well know result, see for
example \cite{Ein}, Proposition 1.1.

\begin{prop}\label{prop: cone iff degenerate}
The hypersurface $V(f)=X$ is a cone if and only if $X^*$ is a
degenerate variety. In particular the hypersurface $X=V(f)$ is a
cone if and only if the partial derivative of $f$ are linearly
dependent.
\end{prop}

Now we recall a problem considered twice by Hesse in \cite{Hesse1}
and \cite{Hesse2}, giving an equivalent geometric formulation of
it. Note that obviously, when $X=V(f)\subset\PP^n$ is a cone, i.e.
up to a linear change of variables $f$ does not depend on all the
variables, then $h_f\equiv 0$. One can ask if the converse
holds.\\

\textbf{Hesse's problem}:  {\it \bf If $h_f\equiv 0$,  then
is $V(f)\subset\PP^n$
a cone?}\\

Note that by Proposition \ref{prop: cone iff degenerate}
Hesse' s claim is equivalent to prove that if $h_f\equiv 0$ then the derivatives of $f$ are linearly dependent.\\
The question was reconsidered by  Gordan and Noether in
\cite{Gordan Noether}. They showed that Hesse's claim is true when
$n\leq 3$ but it is false in general when $n\geq 4$. They
constructed families of counterexamples which have been revisited
recently by Permutti in \cite{Permutti1}, \cite{Permutti2} and by
Lossen in \cite{Lossen}. An easy  example for $n=4$ is the following cubic polynomial
$f(x_0,x_1,x_2,x_3,x_4)=x_0x_3^2+2x_1x_3x_4+x_2x_4^2$.

\begin{rem} {\rm Note that if  $d={\rm deg}(f)\leq 2$ then the Hesse's
claim is true for every $n\geq 1$. Indeed if $d=1$ then $V(f)$ is
a hyperplane, and so it is a cone. If $d=2$ then $V(f)$ is a
hyperquadric and $\ch_f$ is the matrix associated to the quadratic
form of $V(f)$. Since its determinant is zero, the associated
hyperquadric is singular, and so the hyperquadric is a cone. }
\end{rem}

\textbf{From now on $f\in
k[x_0,\ldots,x_n]$ is a homogeneous reduced  polynomial of degree $d\geq 3$ such that
$h_f\equiv 0$.}\\

Since $h_f\equiv 0$, there exist homogeneous polynomials $\pi\in
k[y_0,\ldots , y_n]$ such that $\pi(f_0,\ldots,f_n)\in
k[x_0,\ldots,x_n]$ is identically equal to zero. Let $g\in
k[y_0,\ldots , y_n]$ be such a homogeneous polynomial with this
property and such that $g_i:=\frac{\partial g}{\partial
y_i}(\frac{\partial f}{\partial x_0},\ldots ,\frac{\partial
f}{\partial x_n})\in k[x_0,\ldots , x_n]$, $i=0,\ldots ,n$ are not
all identically equal to zero.

\begin{defi}{\em Let $S=V(g)\subset\PP^{n*}$ be an irreducible and reduced hypersurface  containing the polar image $Z(f)$ and such that $Z(f)$ is not completely contained in the singular locus of $S$.
Let
$$\psi_g\colon\mathbb{P}^n\dashrightarrow\mathbb{P}^n$$
be the composition of $\phi_f$ with $\phi_{g}$ (or equivalently $\psi_g$ is the composition of $\phi_f$ with the Gauss map of $S$).
 If the polynomials
$g_i:=\frac{\partial g}{\partial y_i}(\frac{\partial f}{\partial x_0},\ldots ,\frac{\partial f}{\partial x_n})\in k[x_0,\ldots ,
x_n]$ have a common
divisor $\rho:={\rm g.c.d.}(g_0,\ldots ,g_n)\in k[x_0,\ldots ,
x_n]$,
set $h_i:=\frac{g_i}{\rho}\in k[x_0,\ldots ,
x_n]$, for $i=0,\ldots ,n$.

It follows that the map $\psi_g$ is given by:
$$\psi_g(p) =(g_0(f_0(p),\ldots ,f_n(p))\colon\ldots\colon g_n(f_0(p),\ldots ,f_n(p)))=(h_0(p)\colon\ldots\colon h_n(p)),$$
with ${\rm
g.c.d.}(h_0,\ldots, h_n)=1$.}
\end{defi}

So we have:
$$
\begin{array}{cclcl}
X&\stackrel{\mathcal{G}_X}{\dashrightarrow}&X^*&&\\
\cap&&\cap&&\\
\mathbb{P}^n&\stackrel{\phi_f}{\dashrightarrow}&Z(f)\subset S=V(g)\subset
\mathbb{P}^{n*}&\stackrel{\phi_{g}}{\dashrightarrow}&\mathbb{P}^{n**}\cong \mathbb{P}^n\\ \\
\mathbb{P}^n&--&\stackrel{\psi_g=\phi_g\circ\phi_{f}}{------------}&\dashrightarrow&\mathbb{P}^{n**}\cong
\mathbb{P}^n.
\end{array}
$$
\medskip

Note that the base locus of $\psi_g$ is the scheme
$Bs(\psi_g)=V(h_0,\ldots, h_n)\subset\mathbb{P}^n$ of codimension
at least 2 (because ${\rm g.c.d.}(h_0,\ldots, h_n)=1$ and because
we can assume that the $h_i$'s are not constant).

Set $S^*_Z:=\overline{\psi_g(\mathbb{P}^n)}$ and note that by definition of $\psi_g$,
$$S^*_Z\subset Z(f)^*.$$
Indeed $Z(f)^*$ is made up of the hyperplanes containing the
tangent spaces to $Z(f)$, and $S_Z^*$ is made up by the
hyperplanes which are tangent to $S$ in the points of $Z(f)$.
Since $Z(f)\subset S$ the hyperplanes which are tangent to $S$ in
points of $Z(f)$ are clearly hyperplanes containing the tangent
spaces to $Z(f)$.

Recalling Proposition \ref{prop: dimension of Z dual}, we get:
\begin{equation}\label{image}
S^*_Z\subset Z(f)^*\subset {\rm Sing}(X).
\end{equation}

\medskip
Let us recall a fundamental result proved by Gordan and Noether (see \cite{Gordan Noether} and also \cite{Lossen}, 2.7).

\begin{theorem}\label{theorem: partial diff }
Under the above notation, let $F\in k[x_0,\ldots ,x_n]$. Then:
\begin{equation}\label{2.7}
\sum_{i=0}^n\frac{\partial F}{\partial x_i}h_i=0\,
\Leftrightarrow\, \forall\lambda\in k, \, F(\underline{x})=F(\underline{x}+\lambda\psi_g(\underline{x})).
\end{equation}
\end{theorem}

\begin{rem}\label{rem: f2/xixjhi=0}{\em
Note that $\sum_{i=0}^n\frac{\partial ^2f}{\partial x_i\partial
x_j}h_i=0$. This relation is obtained differentiating the equation
$g(f_0,\ldots,f_n)=0$ with respect to $x_j$ and applying the chain
rule. As a consequence, we get the following relation by Theorem
\ref{theorem: partial diff }:
\begin{equation}\label{2.5}
   f_i(\underline{x})=f_i(\underline{x}+\lambda\psi_g(\underline{x})).
\end{equation}
}\end{rem}

\begin{rem}\label{rem2.8}{\em Using the above result one can find that:
$\sum_{i=0}^n\frac{\partial g_k}{\partial x_i}h_i=0$.\\
Indeed since $g_k$ is a polynomial in $\frac{\partial f}{\partial
x_j}$, $j=0,\ldots, n$,
$$\sum_{i=0}^n\frac{\partial g_k}{\partial x_i}h_i=\sum_{i=0}^n\left(\sum_{j=0}^n\frac{\partial g_k}{\partial (\frac{\partial f}{\partial x_j})}\cdot \left(\frac{\partial^2f}{\partial x_j \partial x_i}\right)\right)h_i=0,$$
where the last equality follows from Remark \ref{rem:
f2/xixjhi=0}.\\
As a consequence we obtain
\begin{eqnarray}\label{eq: hk/dxi=0}\sum_{i=0}^n\frac{\partial h_k}{\partial
x_i}h_i=\frac{1}{\rho}\sum_{i=0}^n\frac{\partial g_k}{\partial
x_i}h_i=0.\end{eqnarray} Since $\psi_g=(h_0:\ldots: h_n)$, by
Theorem \ref{theorem: partial diff } we have
\begin{equation}\label{2.8}
\forall p\in\mathbb{P}^n,\,  \forall\lambda\in k,\,  \psi_g(p)=\psi_g(p+\lambda\psi_g(p)).
\end{equation}}\end{rem}
Geometrically, the equation \eqref{2.8} means that the fiber of
$\psi_g$ over a point $q\in S^*_Z$, $\psi_g ^{-1}(q)$, is a cone
whose vertex contains the point $q$. Indeed $\forall
p\in\mathbb{P}^n$ such that $\psi_g (p)=q$,
$q=\psi_g(p)=\psi_g(p+\lambda \psi_g(p))=\psi_g(p+\lambda q)$,
i.e. $p+\lambda q\in\psi_g^{-1}(q)$ for all $\lambda$. Hence
$\forall p\in\mathbb{P}^n$ such that $\psi_g (p)=q$, the line
$<p,q>$ is contained in $\psi_g ^{-1}(q)$ and $<p,q>\cap {\rm
Bs}(\psi_g)=\{q\}$ as sets.

\begin{rem}\label{rem3}{\em
If the condition \eqref{2.7} holds for a polynomial $F\in
k[x_0,\ldots ,x_n]$, then:
$$S_Z^*\subset V(F).$$
Indeed, using the equation \eqref{2.7} and applying Taylor's formula we have:
$$0= F(\underline{x})-F(\underline{x}+\lambda\psi_g(\underline{x}))=\sum_{k=1}^e\lambda ^k\Phi_k,$$
with $e={\rm deg}(F)$, $\Phi_k:=\sum_{i_1\ldots i_k}\frac{\partial
^kF}{\partial x_{i_1}\ldots \partial x_{i_k}}\frac{h_{i_1}\ldots
h_{i_k}}{k!}.$ In particular, if we assume $F\neq 0$, homogeneous
of degree $e\geq 1$, comparing the coefficient for $\lambda ^e$ we
get: $F(\psi_g (\underline{x}))=F(h_0,\ldots,h_n)=0$. }\end{rem}

Collecting the above remarks we get  that the following result.
\begin{prop}\label{prop: GN dimension of S dual}
{Let notation and hypothesis be as above and suppose that $V(f)$ is not a cone. Then}:
$$S^*_Z \subset V(h_0,\ldots ,h_n)={\rm Bs}(\psi_g).$$
In particular $\dim (S^*_Z)\leq \dim V(h_0,\ldots,h_n)\leq n-2$.
\end{prop}
\bprf By the equation \eqref{eq: hk/dxi=0},
$\sum_{i=0}^n\frac{\partial h_k}{\partial x_i}h_i=0$, hence the
condition \eqref{2.7} holds for $\psi_g=(h_0,\ldots. h_n)$. By
Remark \ref{rem3} this implies that $S_Z^*\subset V(h_0,\ldots,
h_n)=Bs(\psi_g)$.

The bound on the dimension of $V(h_0,\ldots, h_n)$ follows from
the fact that ${\rm g.c.d.}(h_0, \ldots, h_n)=1$ and from the fact
that we can suppose that the $h_i$'s are non-zero and non-costant
since $V(f)$ is not a cone.\eprf

We have the following useful proposition that will be used later.
\begin{prop}\label{prop: GN lines in base locus}
Under the above notation, let $q\in S^*_Z$ be a general point and
let $w\in{\rm Bs}(\psi_g)$ (respectively $t\in {\rm Sing}(X)$). If
$w\in \overline{\psi_g^{-1}(q)}\setminus\{q\}$, (respectively
$t\in \overline{\psi_g^{-1}(q)}\setminus\{q\}$), then the line
$\langle w,q\rangle$ is contained in ${\rm Bs}(\psi_g)$,
(respectively the line $\langle t,q\rangle$ is contained in ${\rm
Sing}(X)$).
\end{prop}
\bprf Since $\psi_g ^{-1}(q)$ is a cone whose vertex contains the point $q$ by \eqref{2.8}, then $\overline{\psi_g^{-1}(q)}$
is a cone whose vertex contains the point $q$. The line $\langle w,q\rangle$, respectively the line  $\langle t,q\rangle$, is contained in $ \overline{\psi_g^{-1}(q)}$, whence the conclusion follows from the relations \eqref{2.8} and \eqref{2.5}.
\eprf

Another general and useful remark is the following lemma which gives a connection between the polar map of the restriction to a hyperplane with the geometry of $Z(f)$ (see \cite{CilibertoRussoSimis}, Lemma 3.10).
\begin{lemma}\label{projection}
Let $X=V(f)\subset\mathbb{P}^n$ be a hypersurface. Let $H=\mathbb{P}^{n-1}$ be a hyperplane not contained in $X$ and let $h=H^*$ be the corresponding point in $\mathbb{P}^{n*}$ and let $\pi_h$ denote the projection from the point $h$. Then:
$$\phi_{V(f)\cap H}=\pi_h\circ (\phi_{V(f)}|H).$$
In particular, $Z(V(f)\cap H)\subset \pi_h(Z(f))$, where  $Z(V(f)\cap H)$ denotes the closure of the image of the polar map $\phi_{V(f)\cap H}$ of the hypersurface $V(f)\cap H$ of $H$.
\end{lemma}
\section{Cases in which Hesse's claim is true.}
In this section we shall consider some hypotheses under which the
conclusion in Hesse' s claim holds.
\begin{rem}\label{rem: S* o Z(f)* cone X cone} {\em
\begin{itemize}
\item [i)] \textit{ Let $S=V(g)\supseteq Z(f)$. If $S^*$ is a
cone, then $X$ is a cone. If $Z(f)^*$ is a
cone then $X$ is a cone.} \\
Indeed if $S^*$ (resp. $Z(f)^*$) is a cone, then $S$ (resp.
$Z(f)$) is a degenerate variety of $\pdual$. Since $X^*\subset
Z(f)\subset S$, $X^*$ is a degenerate variety, whence $X$ is a
cone. \item [ii)] \textit{If $\dim
(S^*)=0$ {\rm (}resp. $\dim(Z(f)^*)=0${\rm)} then $X$ is a cone.}\\
By reflexivity we get $S=S^{**}=\PP^{n-1}$, resp.
$Z(f)=\PP^{n-1}$. In both cases the result follows from part i).\\
Again by part i), \textit{if $S^*$ {\rm(}resp. $Z(f)^*${\rm)} is a linear
subspace of $\mathbb{P}^{n**}$ then $X$ is a cone} (because the
dual of linear subspaces of $\mathbb{P}^{n**}$ are linear
subspaces of $\PP^{n*}$).
\item [iii)] \textit{If $\dim (S^*_Z)=0$, then $X^*$ is a cone.}\\
Indeed if $\dim (S^*_Z)=0$, $S^*_Z$ is a point, and then all the tangent
spaces to the points of $Z(f)$ are contained in a hyperplane (the
dual of the point $S^*_Z$). But this means that $Z(f)$ is contained in
a hyperplane, whence it is degenerate. It follows  that $X^*$ is
degenerate and so $X$ is a cone.
\end{itemize}
}\end{rem}

In particular we recall some properties of the cone $X$ which are
described dually by other geometric  properties of its dual
variety $X^*$.
\begin{rem}\label{rem: X* properties X cone} {\em
\begin{itemize}
\item [i)] \textit{If $X^*$ is a non degenerate subvariety of a
hyperplane $\PP^{n-1}$ in $\PP^{n*}(\cong\PP^{n})$, then $X$ is a
cone with vertex exactly a point.} \item [ii)] \textit{if $X^*$ is
a non-degenerate subvariety of a linear subspace $L=\PP^{n-m}$
($m=1,\ldots n-1$) in $\PP^{n*}(\cong\PP^{n})$, then $X$ is a cone
with vertex a linear subspace $\PP^{m-1}=L^*$.} \item [iii)]
\textit{If $X^*$ is union of $d\geq 1$ points which span a linear
subspace $\PP^{n-m}$ of $(\PP^n)^*$, then $X$ is made up by $d$
hyperplanes whose intersection is a $(m-1)$-linear subspace of
$\PP^{n}$.}
\end{itemize}
}\end{rem}

Now we can prove easily Hesse's claim when $n\leq 3$.

\begin{prop}\label{case1}
Let $X=V(f)\subset\mathbb{P}^1$ be a reduced hypersurface of degree $d$. Then $X=V(f)$ has vanishing Hessian if and only if $X$ is a cone. In this case $d=1$ and $X$ is a point.
\end{prop}
\bprf In this case $Z(f)\subsetneq \mathbb{P}^1$ must be a point
because $\phi_f$ is not dominant, so the partial derivatives of
$f$ are constant and $d=1$ since $X$ is reduced, i.e. $X$ is a
point. \erem

\begin{prop}
Let $X=V(f)\subset\mathbb{P}^2$ be a reduced hypersurface of
degree $d\geq 2$. Then $X=V(f)$ has vanishing Hessian if and only
if $X$ is a cone, i.e. if and only if $X$ consists of $d$ distinct
lines through a point.
\end{prop}
\bprf
Note that $\dim (Z(f))\leq 1$.
As in Proposition \ref{case1}, $Z(f)$ is a point if and only if $d=1$.
Assume $\dim (Z(f))=1$. By Proposition \ref{prop: dimension of Z dual}, $Z(f)^*\subset {\rm Sing}(X)$. Since we are assuming $X$ to be reduced,
we infer that $Z(f)^*$ is a point, so $Z(f)$ is a line and whence the hypersurface $X$ is a cone,  made up by $d$ lines meeting in the point $Z(f)^*$ (where $d$ is the degree of $f$).
\eprf

The following result was proved by Gordan and Noether in
\cite{Gordan Noether}. Here we give an easier and more geometrical
proof of it.
\begin{prop}\label{case n=3}
Let $X=V(f)\subset\mathbb{P}^3$ be a reduced hypersurface of
degree $d\geq 3$. Then $X=V(f)$ has vanishing Hessian if and only
if $X$ is a cone. More precisely, $X=V(f)$ has vanishing Hessian
if and only if either $X$ is a cone over a curve of vertex a point
or $X$ consists of $d$ distinct planes through a line. In the
first case $Z(f)$ is a plane in $\mathbb{P}^{3*}$ while in the
second case it is a line in $\mathbb{P}^{3*}$.
\end{prop}
\bprf In this case $\dim (S^*_Z)\leq\dim (Z(f)^*)\leq 1$  by
\eqref{image} and Proposition \ref{prop: dimension of Z dual}.

Assume $\dim (S^*_Z)=0$ and we are in the hypothesis of the Remark
\ref{rem: S* o Z(f)* cone X cone}  ii), so $X$ is a cone. In
particular if $\dim (Z(f)^*)$=0, $X^*$ is contained in the plane
$Z(f)\subset\mathbb{P}^3$. Moreover $X^*$ is non-degenerate in the
plane $Z(f)$ because otherwise it would be either a line or a
point, which is clearly impossible. It follows from Remark
\ref{rem: X* properties X cone} that $X$ is a cone with vertex the
point $Z(f)^*$ over a plane curve (the dual curve of $X^*$ in the
plane $Z(f)$).

Assume now that $\dim (S^*_Z)=1$. Since $\dim (Z(f)^*)\leq 1$ and
$S^*_Z\subset Z(f)^*$, this implies that $\dim (Z(f)^*)=1$. Since
$Z(f)^*$ and $S^*_Z$ are irreducible (because $Z(f)$ is
irreducible), $Z(f)^*=S^*_Z$.\\ Let ${s_1}$, ${s_2}$ two distinct
general point of $S^*_Z$. Then $\overline{\psi_g^{-1}({s_i})}$ is
a surface which is a cone whose vertex contains the point $s_i$.
Let
${t}\in\overline{\psi_g^{-1}({s_1})}\cap\overline{\psi_g^{-1}({s_2})}\subset
{\rm Bs}(\psi_g)$. By Proposition \ref{prop: GN lines in base
locus},
 the lines $\langle {s_i}, {t}\rangle$,
$i=1,2$, are contained in the base locus of $\psi_g$. Since $\dim
Bs(\psi_g)\leq 1$, the irreducible component of $Bs(\psi_g)$ passing
through ${s_1}$ is exactly the line $\langle
{s_1},{t}\rangle$. But also $S^*_Z$ is an
irreducible component of $Bs(\psi_g)$ of dimension one passing
through ${s_1}$, so it has to coincide with the line
$\langle {s_1}, {t}\rangle$. We conclude that
$S^*_Z=Z(f)^*=\langle {s_i}, {t}\rangle =\langle
{s_1}, {s_2}\rangle$.

Since $Z(f)^*$ is a line, then $Z(f)$ is a line and $X^*\subsetneq
Z(f)=\mathbb{P}^1$, whence $X$ is the union of $d$ planes through
$Z(f)^*=\mathbb{P}^1$ by Remark \ref{rem: X* properties X cone}.
 \erem

\begin{cor}\label{Xcone}
Let $X=V(f)\subset\mathbb{P}^n$, $n\geq 4$ be a reduced hypersurface of degree $d$. If $X=V(f)$ has vanishing Hessian and if ${\rm dim}(Z(f))\leq 2$, then $X=V(f)$ is a cone.
\end{cor}
\bprf Let $H\subset\mathbb{P}^n$ be a general $\mathbb{P}^3$ and
let $h=H^*=\mathbb{P}^{n-4}$. By iterating Lemma \ref{projection}
we deduce that the variety $Z(V(f)\cap H)$ is contained in the
variety   $\pi_h(Z(f))$, whose dimension equals ${\rm dim}(Z(f))$.
Thus  $V(f)\cap H$ has vanishing Hessian because the polar map
$\phi_{V(f)\cap
H}\colon\mathbb{P}^3\dashrightarrow\mathbb{P}^{3*}$ is not
dominant. By Proposition \ref{case n=3} we infer that $V(f)\cap H$
is a cone. By the generality of $H$ we get that
$X=V(f)\subset\mathbb{P}^n$ is a cone. \erem

\section{Gordan--Noether and Franchetta counterexamples to Hesse's conjecture.}
In this section we will describe some examples of hypersurfaces in
$\mathbb{P}^n$, $n\geq 4$, with vanishing Hessian and which are
not cones, following \cite{Gordan Noether} and
\cite[$\S$2]{CilibertoRussoSimis}. Moreover we introduce the
hypersurfaces in $\mathbb{P}^4$ which are counterexamples to
Hesse's claim described by Franchetta (cf. \cite{Franchetta}). We
observe that these hypersurfaces are particular cases of the ones
described by Gordan--Noether.

We briefly recall the results of Gordan--Noether and Permutti in
connection with the Hesse problem, following \cite{CilibertoRussoSimis}.

Assume $n\geq 4$ and fix integers $t\geq m+1$ such that $2\leq
t\leq n-2$ and $1\leq m\leq n-t-1$. Consider forms $h_i(y_0,
\ldots ,y_m)\in k[y_0, \ldots ,y_m]$, $i=0,\ldots ,t$, of the same
degree, and also forms $\psi_j(x_{t+1}, \ldots ,x_n)\in k[x_{t+1},
\ldots ,x_n]$, $j=0,\ldots ,m$, of the  same degree. Introduce the
following homogeneous polynomials all of the same degree:

\begin{displaymath}
Q_\ell(x_0,\ldots ,x_n):={\rm det}
\left(
\begin{array}{ccc}
x_0 & \ldots & x_t \\
\frac{\partial h_0}{\partial \psi_0} & \ldots & \frac{\partial h_t}{\partial \psi_0} \\
\ldots & \ldots & \ldots \\
\frac{\partial h_0}{\partial \psi_m} & \ldots & \frac{\partial h_t}{\partial \psi_m} \\
a_{1,0}^{(\ell)} & \ldots & a_{1,t}^{(\ell)} \\
\ldots & \ldots & \ldots \\
a_{t-m-1,0}^{(\ell)} & \ldots & a_{t-m-1,t}^{(\ell)}
\end{array}
\right)
\end{displaymath}
where $\ell =1,\ldots ,t-m$. Note that $a_{u,v}^{(\ell)}\in k$ for $u=1,\ldots ,t-m-1$, $v=0,\ldots ,t$ and $\frac{\partial h_i}{\partial \psi_j}$ stands for the derivative
$\frac{\partial h_i}{\partial y_j}$ computed at $y_j=\psi_j(x_{t+1}, \ldots ,x_n)$ for $i=0,\ldots ,t$ and $j=0,\ldots ,m$.
Let $s$ denote the common degree of the polynomials $Q_\ell$. Taking Laplace expansion along the first row, one has an expression of the form:
$$Q_\ell =M_{\ell ,0}x_0+\ldots +M_{\ell ,t}x_t,$$
where $M_{\ell ,i}$, $\ell=1,\ldots ,t-m$, $i=0,\ldots ,t$ are homogeneous polynomials  of degree
$s-1$ in $x_{t+1}, \ldots ,x_n$.

Fix an integer  $d>s$ and set $\mu =[\frac{d}{s}]$. Fix biforms $P_k(z_{1}, \ldots ,z_{t-m};x_{t+1}, \ldots ,x_n)$ of bidegree $k, d-ks$, $k=0,\ldots ,\mu$. Finally set
\begin{equation}\label{equGN}
f(x_0,\ldots ,x_n):=\sum_{k=0}^{\mu} P_k(Q_1,\ldots ,Q_{t-m},x_{t+1},\ldots ,x_n),
\end{equation}
a form of degree $d$ in $x_0,\ldots ,x_n$. The polynomial $f$ is called a {\it Gordan--Noether polynomial} (or a {\it GN--polynomial}) of type $(n,t,m,s)$, and so will also any polynomial which can be obtained from it by a projective change of coordinates. Accordingly,
 a {\it Gordan--Noether hypersurface} (or a {\it GN--hypersurface}) of type $(n,t,m,s)$ is the hypersurface $V(f)$, where $f$ is a  GN--polynomial of type $(n,t,m,s)$.

The main point of the Gordan--Noether construction is that a
GN--polynomial has vanishing Hessian. For a proof see
\cite[Proposition 2.9]{CilibertoRussoSimis}. Another proof closer
to Gordan--Noether's original approach is contained in
\cite{Lossen}.

\begin{prop}\label{GNpol}
Every GN--polynomial has vanishing Hessian.
\end{prop}
Following \cite{Permutti2} and \cite{CilibertoRussoSimis} we give
a geometric description of a GN--hypersurface of type $(n,t,m,s)$
as follows. The main result is that the GN--hypersurfaces have
vanishing Hessian (cf. \ref{GNpol}) but in general they are not
cones, so they are counterexample to Hesse's conjecture.

\begin{defi}
{\em Let $f$ be a GN--hypersurface of type $(n,t,m,s)$. The} core {\em of $V(f)$ is the $t$-dimensional subspace $\Pi\subset V(f)$ defined by the equations $x_{t+1}=\ldots =x_n=0$.
}\end{defi}

\medskip
We will call a GN--hypersurface of type $(n,t,m,s)$ {\it general}
if the defining data, namely the polynomials $h_i(y_0, \ldots
,y_m)$, $i=0,\ldots ,t$, the polynomials $\psi_j(x_{t+1}, \ldots
,x_n)$, $j=0,\ldots ,m$ and the constants $a_{u,v}^{(\ell)}$,
$\ell =1,\ldots ,t-m$, $u=1,\ldots ,t-m-1$, $v=0,\ldots ,t$, have
been chosen generically.

\begin{prop}{\rm (\cite[Proposition 2.11]{CilibertoRussoSimis})}
\label{2.11}
Let $V(f)\subset\mathbb{P}^n$  be a GN--hypersurface of type
$(n,t,m,s)$ and degree $d$. Set $\mu =[\frac{d}{s}]$. Then
\begin{itemize}
\item [{\rm i)}] $V(f)$ has multiplicity $d-\mu$ at the general point of its core $\Pi$.
\item [{\rm ii)}] The general $(t+1)$-dimensional subspace $\Pi_\xi\subset\mathbb{P}^n$ through $\Pi$ cuts out on $V(f)$, off $\Pi$, a cone of degree $\mu$ whose vertex is a $m$-dimensional subspace $\Gamma_\xi\subset\Pi$.
\item [{\rm iii)}] As $\Pi_\xi$ varies the corresponding subspace $\Gamma_\xi$ describes the family of tangent spaces to an $m$-dimensional unirational subvariety $S(f)$ of $\Pi$.
\item [{\rm iv)}] If $V(f)$ is general and $\mu > n-t-2$ then $V(f)$ is not a cone.
\item [{\rm v)}] The general GN--hypersurface is irreducible.
\end{itemize}
\end{prop}

\begin{defi}\label{Frhyp}{\em (\cite{Franchetta})
A reduced hypersurface $F=V(f)\subset \mathbb{P}^4$ of degree $d$ is said to be a {\it Franchetta hypersurface} if it is
swept out by a one-dimensional family $\Sigma$ of planes such
that:
\begin{itemize}
\item all the planes of the family $\Sigma$ are tangent to a plane
rational curve $C$  (of degree $p>1$) lying on $F$; \item the
family $\Sigma$ and the curve $C$ are such that for a general
hyperplane $H=\PP^3\subset\PP^4$ passing through $C$, the
intersection $H\cap F$, off the linear span of $C$, is the union
of planes of $\Sigma$ all tangent to the curve $C$ in the same
point $p_H$. \erem\end{itemize} }\end{defi}

\begin{rem}\label{remark: Franchetta equvalent to GN}{\em
Note that by Proposition \ref{2.11}  a GN--hypersurface
$X=V(f)\subset\mathbb{P}^4$ of type $(4,2,1,s)$ is a Franchetta
hypersurface with core the linear span of the curve $C$. On the
contrary Permutti proved in \cite{Permutti1} that  a Franchetta
hypersurface $V(f)\subset\mathbb{P}^4$ is a GN--hypersurface of
type $(4,2,1,s)$. In particular (by Proposition \ref{GNpol}) a
Franchetta hypersurface has vanishing Hessian. This fact can be
proved directly see also \cite{Permutti1} and \cite[Proposition
2.18]{CilibertoRussoSimis}.}\end{rem}

\section{A geometrical proof of Gordan--Noether and Franchetta classification of hypersurfaces in $\mathbb{P}^4$ with vanishing Hessian.}
In the previous section we saw that the classes of
GN-hypersurfaces of type $(4,2,1,s)$ and of Franchetta
hypersurfaces coincide. In this section we use the geometrical
methods developed in the first section and some other easy facts
to provide a short and selfcontained proof of Franchetta
characterization of hypersurfaces with vanishing Hessian in
$\mathbb{P}^4$, \cite{Franchetta}. So we will prove in a
geometrical way that the hypersurfaces in $\mathbb{P}^4$ with
vanishing Hessian are either cones or Franchetta hypersurfaces and
that there are no other possibilities. A similar result is not
known in higher dimension.\smallskip

First we give a preliminary result describing a geometrical
consequence of the vanishing of the hessian of hypersurfaces in
$\PP^4$, not cones.

\begin{prop}\label{case n=4}
Let $X=V(f)\subset\mathbb{P}^4$ be a reduced hypersurface of
degree $d\geq 3$, not a cone. If $X=V(f)$ has vanishing Hessian
then $Z(f)^*\subset\PP^4$ is an irreducible plane rational curve.
 Equivalently
$Z(f)$ is a  cone with vertex a line  over an irreducible plane
rational curve.
\end{prop}
\bprf By Corollary \ref{Xcone}, we can suppose ${\rm
dim}(Z(f))\geq 3$. Thus $Z(f)^*=S^*_Z$, and $Z(f)=V(g)$  with
$g\in k[y_0,\ldots ,y_4]$ an irreducible polynomial. Note  that by
Proposition \ref{prop: GN dimension of S dual} we have $1\leq \dim
(Z(f)^*)\leq \dim({\rm Bs}(\psi_g))\leq 2$.

\smallskip Assume first $\dim (Z(f)^*)=2$ so that $Z(f)^*$ is an irreducible component of
${\rm Bs}(\psi_g)$. Consider the intersection between the closure
of the fibers on two different general points, ${s_1},\ {s_2}\in
Z(f)^*$. The fiber on each of these points has dimension two, so
there exists
${t}\in\overline{\psi_g^{-1}({s_1})}\cap\overline{\psi_g^{-1}({s_2})}$.
By Proposition \ref{prop: GN lines in base locus}, the lines
$\langle {s_i},{t}\rangle$, $i=1,2$, are contained in $Bs(\psi_g)$
and hence in the irreducible component of it containing $s_1$ and
$s_2$. Since $s_1$ and $s_2$ are general points in $Z(f)^*$,
$Z(f)^*$ is the unique irreducible component of $Bs(\psi_g)$
containing them. Furthermore  $Z(f)^*$ is a ruled surface (because
through a general point $s\in Z(f)^*$ there passes a line $\ell_s$
contained in $Z(f)^*$), which is a  cone (because $\ell
_{s_1}\cap\ell_{s_2}\neq\varnothing$ for $s_1,s_2\in Z(f)^*$
general points), whence by Remark \ref{rem: S* o Z(f)* cone X
cone}, $X$ is a cone.\\
Thus we can assume $\dim (Z(f)^*)=1$. Let ${s_1}$ and ${s_2}$ be
two general points of $Z(f)^*$. Then the intersection
$\overline{\psi_g ^{-1}({s_1})}\cap\overline{\psi_g ^{-1}({s_2})}$
is a surface, say $R$, contained in ${\rm Bs}(\psi_g)$. Note that
this intersection has to stabilize for general points of $Z(f)^*$.

Furthermore for every point $t\in R$ and for a general point $s\in
Z(f)^*$, by Proposition \ref{prop: GN lines in base locus}, the
line $\langle {s}, {t}\rangle$ is contained in ${\rm
Bs}(\psi_g)\cap\overline{\psi_g ^{-1}({s})}$, and hence in $R$. It
follows that $Z(f)^*$ is contained in the vertex of the surface
$R$, and that $R$ (and in fact the intersection of two general
fibers of $\psi_g$) is a plane  ($Z(f)^*$ is not a
line by assumption, so it  cannot be contained in the intersection  of two or more planes).\\
In other words  $Z(f)^*$ is a plane curve, whose linear span
$\Pi=R$ is an irreducible component of ${\rm Bs}(\psi_g)$. Since
$S^*_Z=Z(f)^*$, Proposition \ref{prop: GN lines in base locus} and
the same argument used above show that $Z(f)^*$ is contained in
the vertex of ${\rm Sing}(X)$. Thus the irreducible components of
${\rm Sing(X)}$ of dimension 2 are planes containing $Z(f)^*$ so
that there is a unique irreducible component of ${\rm Sing(X)}$
which is
a plane: the linear span of $Z(f)^*$, i.e. $\Pi$.\\
Note also that $Z(f)^*$  is rational. In fact the map $\psi_g$ is
a rational dominant map from $\mathbb{P}^4$ to $Z(f)^*$, so
$Z(f)^*$ is a unirational curve and hence a rational curve.

Since $Z(f)^*=S^*_Z\subset \Pi =\mathbb{P}^2$ is an irreducible
rational plane curve (not a line),  $Z(f)$ is a cone of vertex a
line $L=\Pi^*=\mathbb{P}^1$ over an irreducible plane curve
$\Gamma$ (of degree $\geq 2$), which is the dual curve of $Z(f)^*$
in the plane $\Pi$. Furthermore $\Gamma$ is  a rational curve
because in this case the Gauss map of the curve ${Z(f)^*}$ is
birational. \eprf

The description given in Proposition \ref{case n=4} is crucial to
prove that a projective hypersurface $X=V(f)$ in $\mathbb{P}^4$
with vanishing Hessian which is not a cone is a Francehtta
hypersurface. The following result finally gives a
characterization of hypersurfaces in $\mathbb{P}^4$ with vanishing
Hessian, which are not cones.

\begin{theorem}\label{HGNF}
Let $X=V(f)\subset\mathbb{P}^4$ be an irreducible and reduced hypersurface of degree $d\geq 3$, not a cone.
The following conditions are equivalent:
\begin{itemize}
\item [{\rm i)}] $X=V(f)$ has vanishing Hessian. \item [{\rm ii)}]
$X=V(f)$ is a Franchetta hypersurface. \item [{\rm iii)}]
$X^*=V(f)^*$ is a scroll surface of degree $d$, having a line
directrix $L$ of multiplicity $e$, sitting in a $3$-dimensional
rational cone $W(f)$ with vertex $L$, and the general plane ruling
of the cone cuts $V(f)^*$ off $L$ along $\mu \leq e$ lines of the
scroll,  all passing through the same point of $L$. \item [{\rm
iv)}]  $X=V(f)$  is a general GN--hypersurface of type
$(4,2,1,s)$, with $\mu =[\frac{d}{s}]$, which has a plane of
multiplicity $d-\mu$.
\end{itemize}
In particular, $X^*=V(f)^*$ is smooth if and only if  $d=3$, $X^*=V(f)^*$ is a rational normal scroll of degree $3$ and $X=V(f)$ contains a plane, the orthogonal of the line directrix of $X^*=V(f)^*$, with multiplicity $2$.
\end{theorem}

\bprf Note that  conditions ii) and iii) are easily seen to be
equivalent (the directrix line $L$ of $X^*$ is the dual of the
plane which is the linear span of the curve $C$ of the Franchetta
hypersurface). By Remark \ref{remark: Franchetta equvalent to GN},
the equivalence of  ii) and iv) is clear. The conditions iv)
implies the condition i) by Proposition \ref{GNpol}.  Thus to finish  the
proof it is sufficient to prove that a hypersurface
$X=V(f)\subset\mathbb{P}^4$ with vanishing Hessian, not a cone,
is  a Franchetta hypersurface.

By Proposition \ref{case n=4}, we have that $Z(f)^*\subset {\rm
Sing}(X)\subset X=V(f)$ is an irreducible plane rational curve,
whose linear span is a plane $\Pi =\mathbb{P}^2$. Equivalently,
$Z(f)$ is a cone of vertex the line $L=\Pi^*=\mathbb{P}^1$ over an
irreducible plane curve $\Gamma$, the dual of $Z(f)^*$ as a plane
curve.

Consider now   a general hyperplane $H\subset\mathbb{P}^4$
passing through  the plane $\Pi$ (and not contained in $X=V(f)$).
The intersection $X\cap H$ is a hypersurface in $H=\mathbb{P}^3$
containing the plane $\Pi$ with a certain multiplicity $\mu\geq 0$
and reduced elsewhere. Note also that the point $h=H^*\in L=\Pi^*$
(because $\Pi\subset H$), whence $\pi_h(Z(f))$ is a surface
naturally embedded  in the dual space of $H$. More precisely
 $\pi_h(Z(f))$ is a cone  with vertex  the point
$p_L=\pi_h(L)$ over the plane  curve $\hat{\Gamma}=\pi_h(\Gamma)$.

By Lemma \ref{projection} we infer that $Z(V(f)\cap H)\subset
\pi_h(Z(f))\subset\mathbb{P}^{3*}$, whence (see Proposition
\ref{prop: h=0 iff non dominant}) the hypersurface $V(f)\cap
H\subset H=\mathbb{P}^3$ has vanishing Hessian. By Proposition
\ref{case n=3} it follows that $V(f)\cap H$ is a cone and either
it is a cone over a curve of vertex a point or $V(f)\cap H$
consists of distinct planes passing through a line. In the first
case $Z(V(f)\cap H)$ is a plane in $\mathbb{P}^{3*}$ but this is
not possible because the cone $\pi_h(Z(f))$ is non degenerate.
Therefore $Z(V(f)\cap H)$ is a line in $H^*$ and $V(f)\cap H$ is a
union of planes through the line $T=Z(V(f)\cap H)^*\subset H$,
where duality is considered between $H$ and $H^*$. Since the
hyperplane section $V(f)\cap H$ is singular and since $H$ was
general through $\Pi$, we deduce that  $L=\Pi^*\subset X^*$.

Note that, by  Lemma \ref{projection}, $Z(V(f)\cap
H)=\phi_{V(f)\cap H}(H)=\pi_h(\phi_f(H))$ is a line contained in
$\pi_h(Z(f))$, whence $\phi_f(H)$ is a plane of the ruling of
$Z(f)$ corresponding to a point $y\in\Gamma$. Furthermore the
lines $L_j:=\Pi_j^*$, duals to the planes in $V(f)\cap H$
different from $\Pi$, pass all through the point $h=H^*$ and are
contained in the plane $T^*=\phi_f(H)$ and in $X^*$.

Let $z=\psi_g(H)\in Z(f)^*$. Then $\phi_f(H)^*=T_z(Z(f)^*)=T$,
i.e. the line of intersection of the planes in $V(f)\cap H$  is
the tangent line to the plane curve $Z(f)^*$ in the point $z$. In
conclusion $X=V(f)\subset\mathbb{P}^4$ is a Franchetta
hypersurface, where we can take as the one dimensional family
$\Sigma$ of planes contained in $X$ exactly the intersection of a
general $\mathbb{P}^3$ through $\Pi$ with $X=V(f)$ (i.e.  the
intersection of the fibers of $\psi_g$ with $X=V(f)$) and we consider as the curve
$C$ (cf. Definition \ref{Frhyp}) the curve $Z(f)^*$.\erem

\addcontentsline{toc}{section}{ \hspace{0.5ex} References}

\end{document}